# PARAMETER SELECTION IN PARTICLE SWARM OPTIMIZATION FOR TRANSPORTATION NETWORK DESIGN PROBLEM

Mehran Fasihozamn Langerudi
Ph.D. Candidate
Department of Civil and Materials Engineering
University of Illinois at Chicago
842 W. Taylor St.
Chicago, IL 60607
Phone: 312-996-0962
Fax: 312-996-2426
Email: mfasih2@uic.edu

In transportation planning and development, transport network design problem seeks to optimize specific objectives (e.g. total travel time) through choosing among a given set of projects while keeping consumption of resources (e.g. budget) within their limits. Due to the numerous cases of choosing projects, solving such a problem is very difficult and time-consuming. Based on particle swarm optimization (PSO) technique, a heuristic solution algorithm for the bi-level problem is designed. This paper evaluates the algorithm performance in the response of changing certain basic PSO parameters.

***Keywords:*** *Transportation; Network Design; optimization; Particle Swarm; inertia weight; constant values*

## 1. Introduction

Transportation Network Design is the important issue of improving transportation networks by selecting the optimal projects among a set of alternatives. TNDP attempts to optimize certain objectives under resource constraints. For an n-project case, considering an accept-reject decision for each project, there are $2^n$ alternative networks which are to be compared. Although solving such a problem among a few alternatives does not take too much time, the solution becomes excessively time-consuming as n increases. Various approaches have been proposed to solve TNDP. In large scale problems, meta-heuristic techniques which benefit some sort of intelligence in finding the optimal solution have proved to be efficient. Particle Swarm

Optimization (PSO) as one of these techniques has already been adapted to TNDP. In this paper, after adapting the PSO to TNDP on the well-known Sioux Falls network, the role of basic PSO parameters in algorithm performance is experimented and the results are shown consequently.

## 2. The TNDP

Let $G = (V, A)$ be a graph representing a transportation network with node set $V$ and arc set $A$, and define $P \subseteq \{(r, s) \in V \times V : r \neq s\}$ as the set of origin-destination (OD) pairs. For each OD pair $(r, s) \in P$, there is a nonnegative flow rate (travel demand) from $r$ to $s$, denoted by $d_{rs}$. In order to simplify the presentation, suppose that $G$ is strongly connected, that is each node $j$ can be reached from every other node $i$ by following a directed path in $G$, and let $K_{rs}$ be the non-empty set of paths from the origin $r$ to the destination $s$.

Define $\bar{A}$ ($\bar{A} \neq A$) as the set of project arcs, and let the decision vector be $y = (y_a)_{a \in \bar{A}}$ with $y_a$ being the binary project decision variable, taking values 0 or 1 depending on rejection or acceptance of any project $a \in \bar{A}$. For a given vector $y$, define the decision network $G_{(y)} = (V, A(y))$ with $A(y) = A \cup \{a \in \bar{A} : y_a = 1\}$ as the set of arcs followed by decision $y$, and for each $(r, s) \in P$ denote by $K_{rs}(y)$ the set of paths joining $r$ to $s$ in $G_{(y)}$. For each path $k \in K_{rs}(y)$ let $f_k$ be the flow of path $k$ from origin $r$ to destination $s$. Moreover, let $\delta_{ak}$ equals 1 if arc $a \in A(y)$ lies on path $k$, and 0 otherwise.

Assume further that each arc $a \in A \cup \bar{A}$ has a node creasing and continuously differentiable travel time function $t_a(x_a): [0, \infty) \to [0, \infty)$ with $x_a$ being the flow rate assigned to arc $a$. Then, letting $c_a$ be the construction cost of project arc $a \in \bar{A}$, and considering the total construction cost being limited to the level of Budget $B$, the TNDP can be illustrated with upper level problem, ULP:

[ULP] $\quad \underset{y}{Min} \ T(y) = \sum_{a \in A(y)} x_a t_a(x_a)$

s.t. $\quad \sum_{a \in \bar{A}} x_a t_a(x_a)$

$\quad\quad y_a = 0$ or $1 \quad \forall a \in \bar{A}$

$\quad\quad X(y)$ is a solution of [LLP(y)]

Where $x(y) = (x_a)_{a \in A(y)}$ is the user equilibrium flow in the decision network $G(y)$, given as the solution of the lower level (traffic assignment) problem, LLP(y), for given $y$:

[LLP(y)] $\quad Min \ \sum_{a \in A(y)} \int_0^{x_a} t_a(w) dw$

s.t. $\quad \sum_{k \in k_{rs}(y)} f_k = d_{rs} \quad \forall (r,s) \in P$

$\quad\quad f_k \geq 0 \quad \forall k \in K_{rs}(y), \forall (r,s) \in P$

$$x_a = \sum_{(r,s)\in P} \sum_{k\in K_{rs}} f_k \delta_{ak} \quad \forall a \in A(y)$$

This is a well-known bi-level programming problem, where the [ULP] seeks a decision vector $y$ for minimizing the total travel time $T(y)$ of the (assigned) traveler, and the [LLP(y)] is the traffic assignment model which estimates the traveler flows, given the decision $y$.

## 3. Particle Swarm Optimization

Particle Swarm Optimization (PSO) is a meta-heuristic optimization approach which has been widely applied to various problems. PSO technique that was developed by Kennedy and Eberhart is originated from the behavior of birds 'flocks in which individuals convey information between themselves and the leader in order to seek the best direction to food.

In a problem space, each particle has a position and a velocity and it moves in the search space with the velocity according to its own previous best position and the group's previous best position. The dimension of the search space can be any positive number. Considering D as the dimension of the search space, the $i^{th}$ particle's position and velocity are represented as $P_i = (p_{ij})_{j=1,\ldots,D}$ and $V_i = (v_{ij})_{j=1,\ldots,D}$ respectively. Each particle maintains its own best position so far achieved as $P^*_i = (p^*_{ij})_{j=1,\ldots,D}$ and the global best position so far recorded by the population as $P^*_g = (p^*_{gj})_{j=1,\ldots,D}$.

During the iteration time t, the velocity of the $j^{th}$ dimension of each particle i is updated by:

$$v_{ij}(t+1) = wv_{ij}(t) + c_1 r_1 (p^*_{ij}(t) - p_{ij}(t)) + c_2 r_2 (p^*_{gj}(t) - p_{ij}(t))$$

Where $w$ is called as the inertia weight, $c_1$ and $c_2$ are constant values and $r_1$, $r_2$ are random numbers in the interval $[0,1]$. The current position of each particle is then defined by the sum of its current velocity and its previous position.

$$p_{ij}(t+1) = p_{ij}(t) + v_{ij}(t+1)$$

In order to avoid the particles from moving out of the search space, the maximum velocity during the iterations is restricted by $v_{max}$.

## 4. Adapting the PSO to the TNDP

Employing the PSO for solving TNDP needs some modifications to the algorithm given in the previous section. First, the PSO is basically developed for continuous optimization problems. This is while the TNDP is formulated as a combinatorial optimization problem in terms of variables $y$ denoted as $|\overline{A}|$-bit binary strings. To adapt the algorithm for this combinatorial nature, one may provide some mapping from the one-dimensional real-valued space to the $|\overline{A}|$-dimensional binary space. This is done here by transforming each real number $p_i$ to its nearest integer in $\left[0, 2^{|\overline{A}|} - 1\right]$, and then transforming the resulting integer in to the base-2 number system as an $|\overline{A}|$-bit binary code. To facilitate the presentation, the latter transformation is illustrated by the function $y(p_i) : \left[0, 2^{|\overline{A}|} - 1\right] \subset Z \rightarrow \{0,1\}^{|\overline{A}|}$.

The PSO must also be adapted for budget constraint embedded in the [ULP].

### 4.1. PSO algorithm

Step 1. Initialization

Select the particle swarm size $n$, the parameters $c_1$ and $c_2$, the value of the inertia weight $w$, and the maximum velocity $v_{max}$.

For $i=1$ to $n$ do: initialize the decision variable $p_i$ so that $\sum_{a \in \bar{A}} c_a y_a(p_i) \leq B$; set $p_i^* = p_i$ and $v_i = 0$.

Set $p_g^* = \arg \min (f(p_1),..., f(p_n))$. Set the iteration counter $t = 0$.

Step 2. Updating each particle's position and velocity

For $i=1$ to $n$ do: generate random numbers $r_1$ and $r_2$ in [0, 1]; update

$v_i \leftarrow w v_i + c_1 r_1 (p_i^* - p_i) + c_2 r_2 (p_g^* - p_i)$; clamp in $v_i$ between the range $[-v_{max}, v_{max}]$ as $v_i = \text{sign}(v_i) \min(|v_i|, v_{max})$; update $p_i \leftarrow p_i + v_i$; transform $p_i$ to its nearest integer in $[0, 2^{|\bar{A}|} - 1]$.

Step 3. Calculating each particle's fitness value

For $i=1$ to $n$ do: set $y = y(p_i)$; if $\sum_{a \in \bar{A}} c_a y_a > B$ then set $f(p_i) = M$ (large fitness value); else, solve the user equilibrium problem [LLP ($y$)] to compute $T(y)$, and set $f(p_i) = T(y)$.

Step 4. Updating local bests and global best

For $i=1$ to $n$ do: update $p_i^* \leftarrow \arg \min (f(p_i^*), f(p_i))$.

Update $p_g^* \leftarrow \arg \min (f(p_g^*), f(p_1),..., f(p_n))$.

Step 5. End criterion.

Set $t = t+1$. If end criterion is not met, go to Step 2. Otherwise, $y = y(p_g^*)$ is the best solution found so far with the objective function value $T(y) = f(p_g^*)$ Collect the necessary information and stop.

## 5. Sioux Falls Network

The Sioux Falls network has 24 nodes and 76 arcs, as shown in Fig. 1. The parameters of the travel time function $t_a(x_a) = \alpha_a + \beta_a x_a^4$ for each arc $a$, and the OD (origin/destination) demands are basically those given in Poorzahedy and Turnquist (1982), and LeBlanc (1975), and are eliminated here for brevity.

There are 10 pairs of project arcs ($|\bar{A}|=10$), of which 5 projects are improvement on existing arcs, and 5 are new arcs. The construction costs of the projects 1-10 are, respectively, 625, 650, 850, 1000, 1200, 1500, 1650, 1800, 1950, and 2100 units of money (Poorzahedy and Abulghasemi 2005). Considering 10 projects, there are $2^{10} (=1024)$ alternative networks. A complete enumeration was used to compute the optimal solution of the TNDP for any given budget level for checking purposes (Poorzahedy and Abulghasemi 2005; Poorzahedy and Rouhani 2007).

## 6. Computational Results

In order to examine the sensitivity of parameters $c_1$, $c_2$ and $w$ in the adapted PSO algorithm, a computer program was implemented in visual basic 6.0 on a laptop with Intel core 2 due 2.4 GHz processor. In this program, the PSO algorithm is terminated after a fixed number of 1000 iterations. Due to the stochastic nature of PSO, the algorithm has been solved 50 times and the results are based on the average values of the 50 runs. As proposed by Hong Zhang, et al 2005, the maximum velocity ($v_{max}$) is set to $x_{max}$. This results in moving more effectively in the search space and accordingly better algorithm performance.

First coefficient $w$ is set as a constant ($=1.1$) and the results are depicted in the space of $c_1$, $c_2$ in figures 2-4. As it can be seen in figure 2, average Number of Traffic Assignment Problems Solved (NTAPS) increases when $c_1$, $c_2$ extends to 2 and then it decreases. In figure 3, the least values of Objective Function Value (OFV) are placed in the area where $c_1 \approx c_2 \approx 2$. Although ascent in average NTAPS where $c_1 \approx c_2 \approx 2$ might imply more algorithm computation time, significant decrease in average OFV speeds up the convergence of the algorithm.

The difference between the average OFV of the first and the last iterations of the algorithm is another value that is shown in figure 4. This value can somehow show the ability of the algorithm to converge since the values of the first iteration are the same for $c_1$, $c_2$ values. Since $c_1$, $c_2$ demonstrate a somewhat symmetric pattern in the figures 2-4, it is reasonable to apply $c_1 = c_2 = c$ in the following discussions. This would decrease the dimension of the discussed parameters and therefore make the following comparisons more sensible.

Figures 5-9 display average NTAPS, average OFV, frequency of finding , difference of OFV between the first and the last iterations and the probability of finding the optimal solution in 50 runs in the space of c, $w$ respectively.

From Figure 5, considering a constant $w$, average NTAPS increases with c until it reaches to its maximum value while average OFV decreases as displayed in Figure 6.

Figure 7 shows the frequency of finding the optimal solution in 1000 iterations that is a determinant of algorithm speed in finding the optimal solution. When $c$ is set in the range $[1.8, 2.2]$, the algorithm reaches the optimal solution in less than 100 iterations. Also, in figure 8, the highest difference between OFV of the first and the last iterations is gained in the same range of $c$. This point can be further approved in figure 9 when this range achieves the highest probability of finding the optimal solution in 50 runs.

From figures 5-9, it can be readily concluded that $w \leq 0.5c + 1$ is a more reliable range for $w$. As a result, by fixing $c = 2$ as the center of the proposed range, figure 10 is drawn to find the best results of $w$ in the range $[0, 2]$. From this figure, it is clearly seen that by selecting $w$ in the range $[0.5, 0.9]$, the decreasing speed of OFV is more considerable than other cases. Referring to figure 9 again, $w = 0.7$ seems to be the best $w$ in this problem.

In many cases, a decreasing variable value for w is proposed. Therefore, a comparison between a constant $w$ (=0.7) and a decreasing $w$ starting from 1.2 to 0.4 is conducted and the results are depicted in figure 11. From this figure, $w$ constant demonstrates a better convergence behavior than $w$ decreasing.

## 7. Summary and Conclusion

By reviewing the papers related to PSO algorithm which were used to solve various problems, It can be figured out that the best PSO coefficients are gained on the basis of specific problem features.

The obtained results of this paper shows that the best parameters for solving TNDP with PSO algorithm are $w = 0.7$, $c_1 = c_2 = 2$. It must be considered that in some cases, the best solution is the one which has the minimum quantity of average NTAPS that is gained by increasing $c_1$, $c_2$.

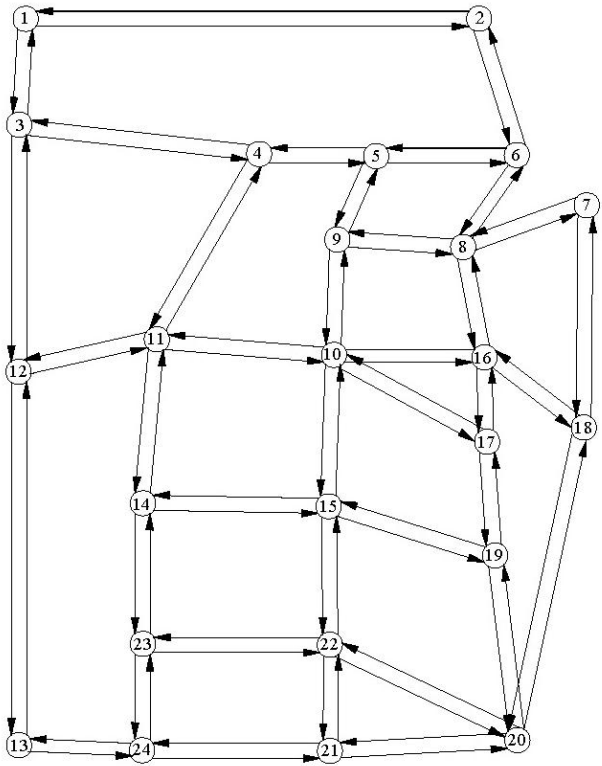

Figure 1: the sioux falls network

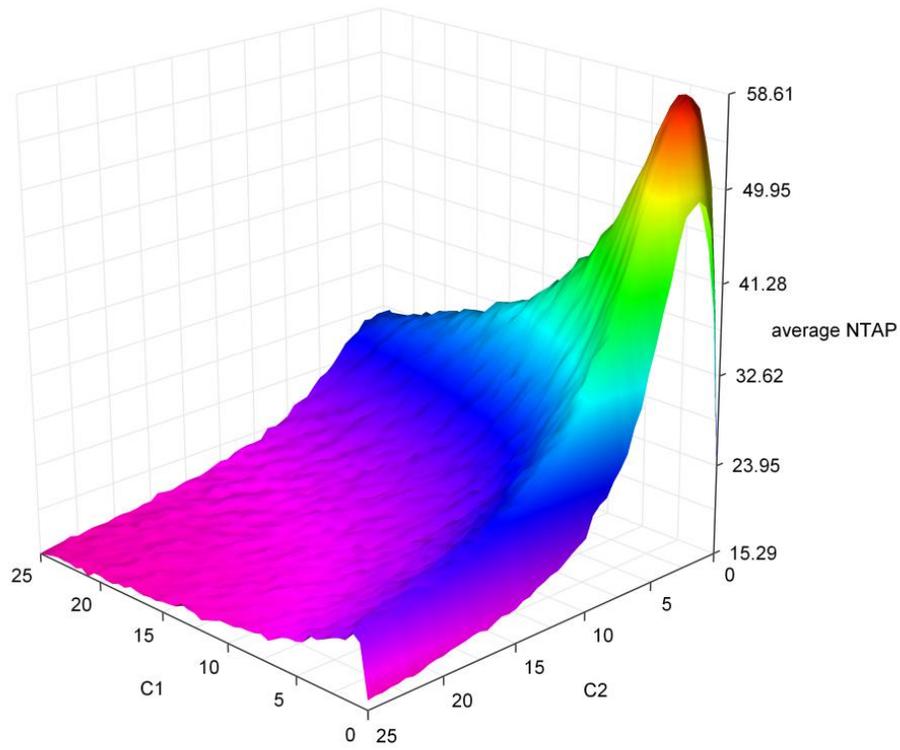

Figure 2: Average Number of Traffic Assignment Problem solved (NTAP)

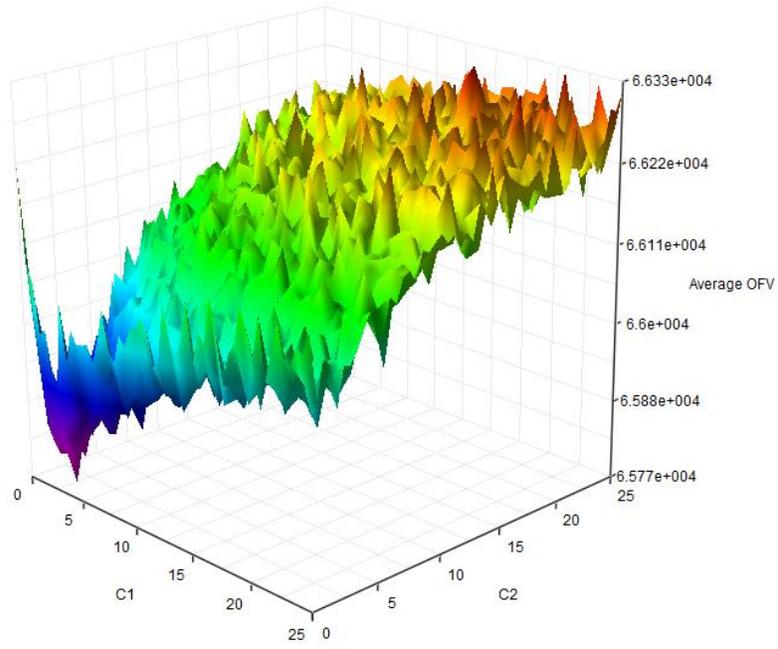

Figure 3: Average Objective Function Value (OFV)

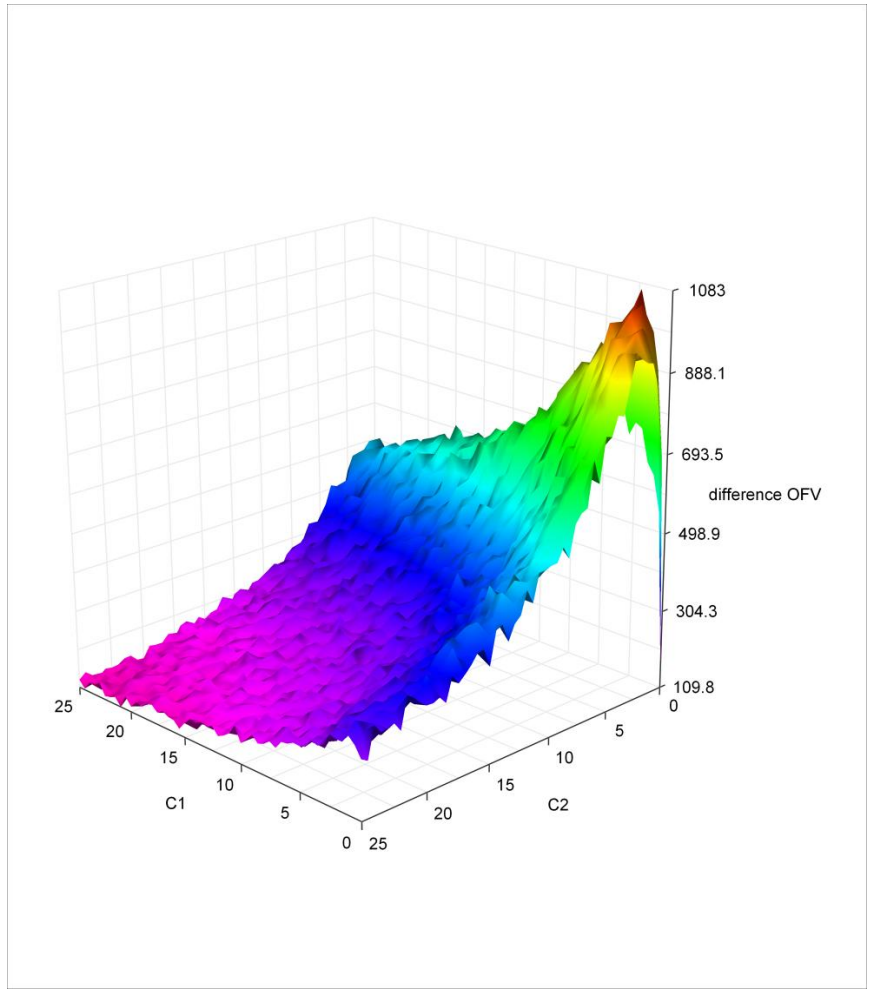

Figure 4: Difference between the average OFV of the first and the last iterations

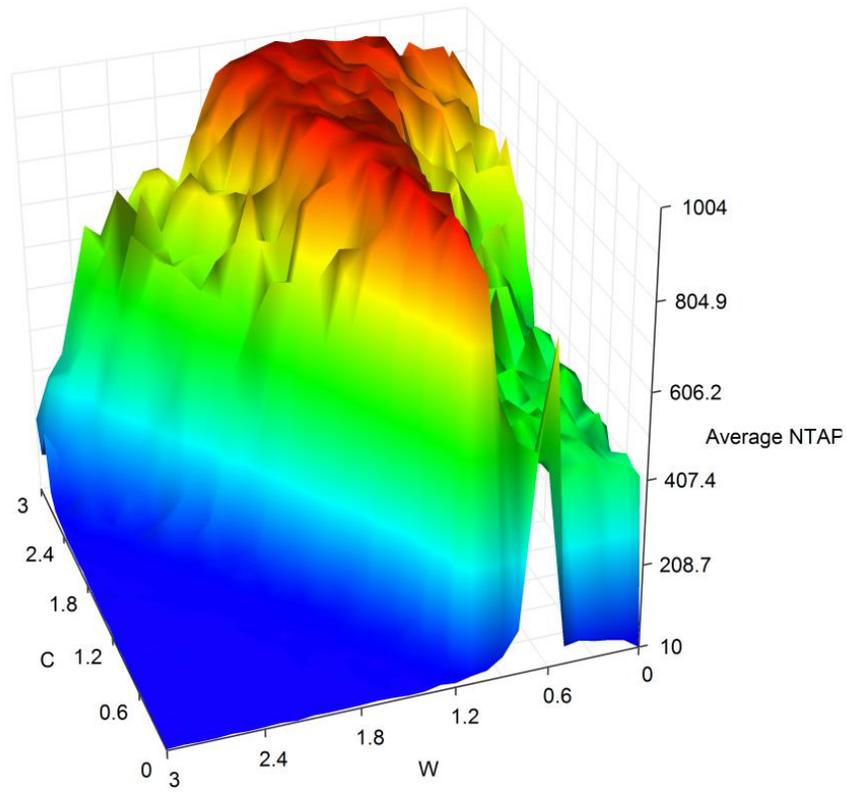

Figure 5: Average Number of Traffic Assignment Problem solved

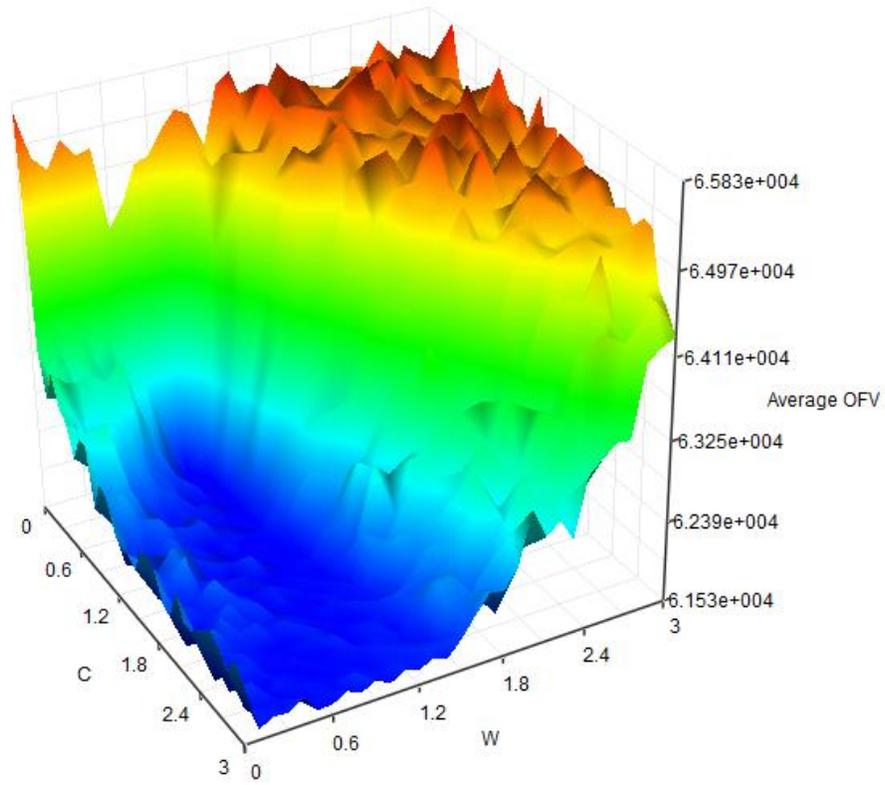

Figure 6: Average Objective Function Value (OFV)

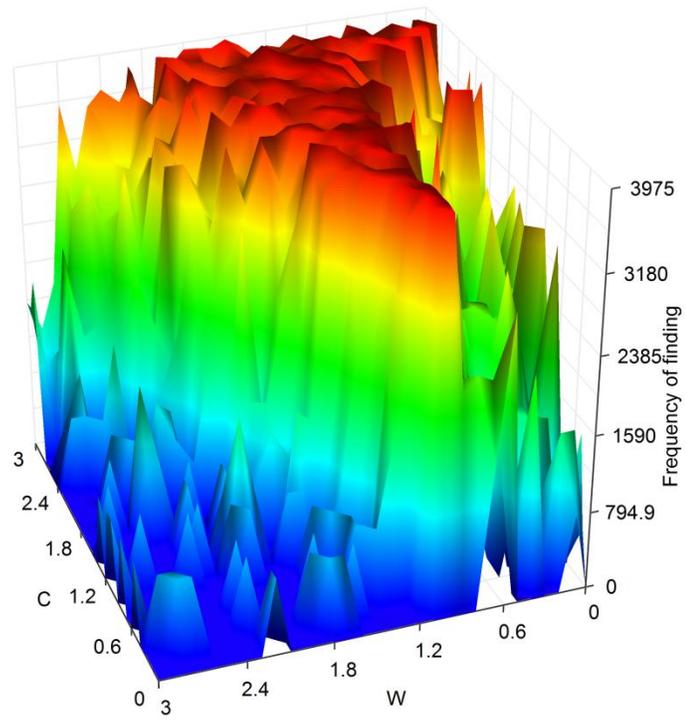

Figure 7: Frequency of finding the optimal solution

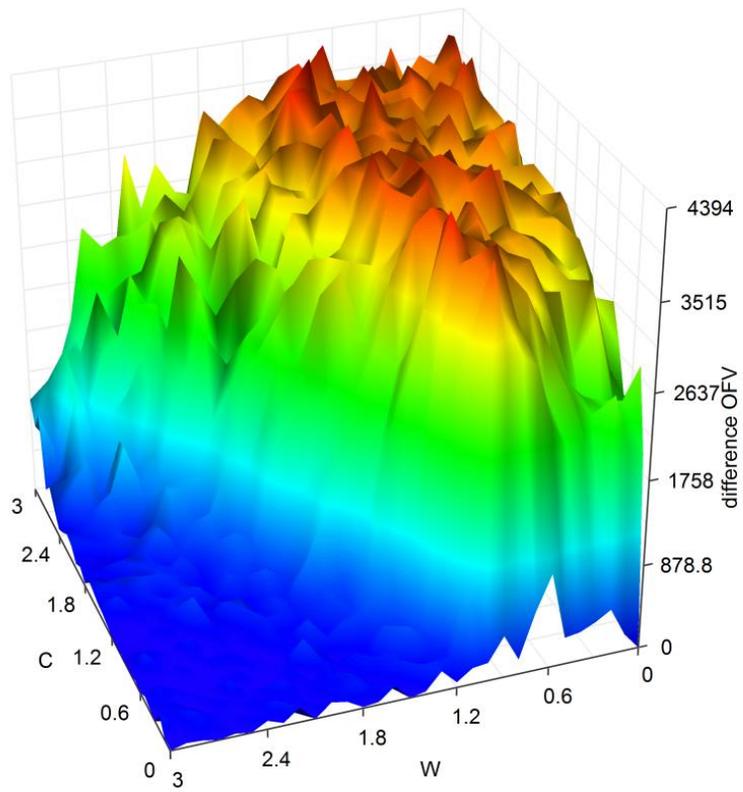

Figure 8: Difference between the average OFV of the first and the last iterations

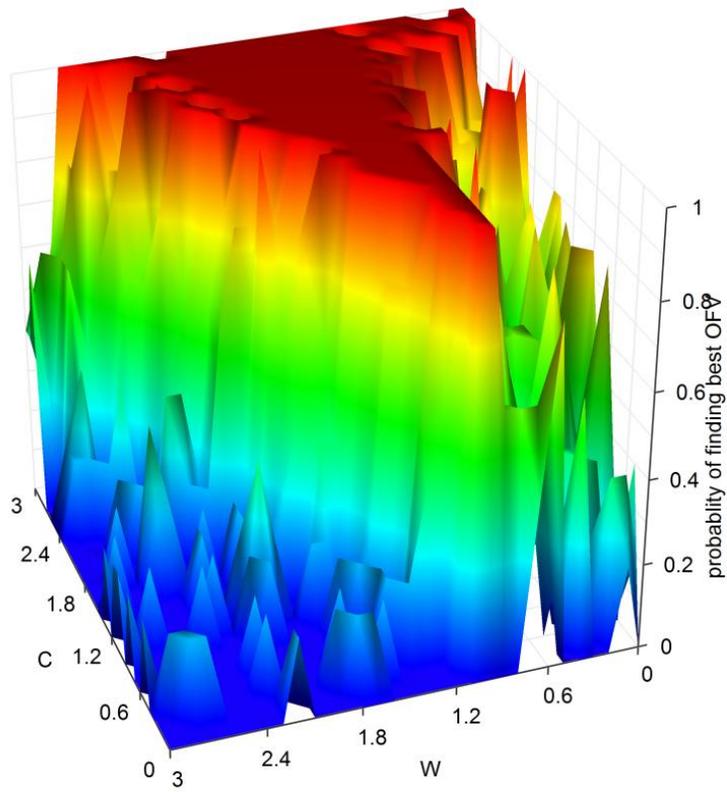

Figure 9: Probability of finding the optimal solution

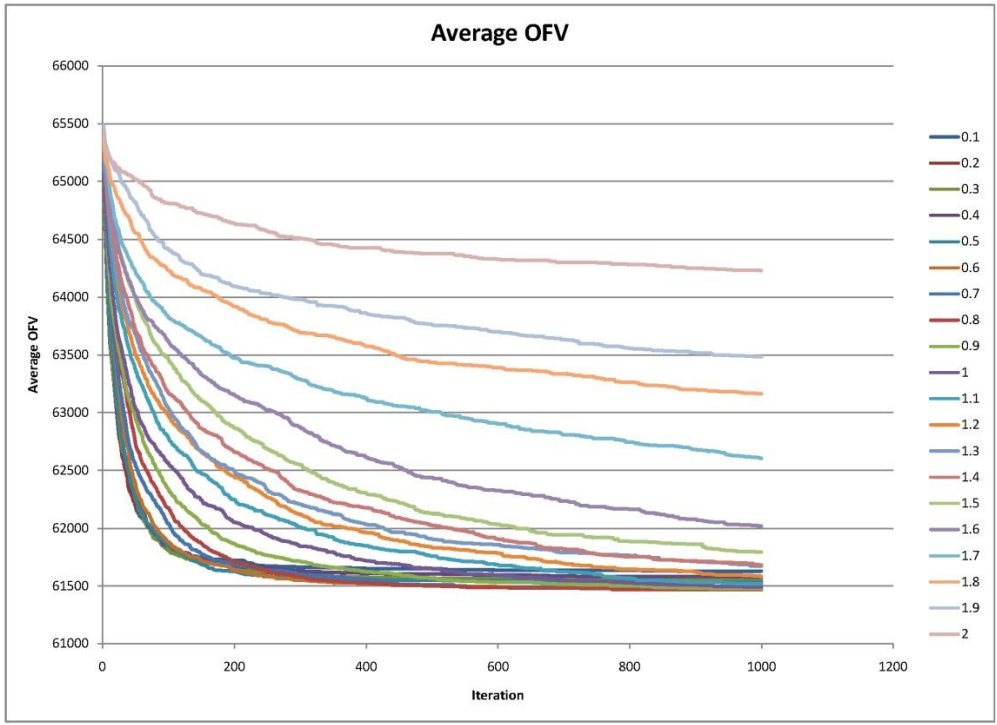

Figure 10: Average Objective Function Value for w in the range [0, 2]

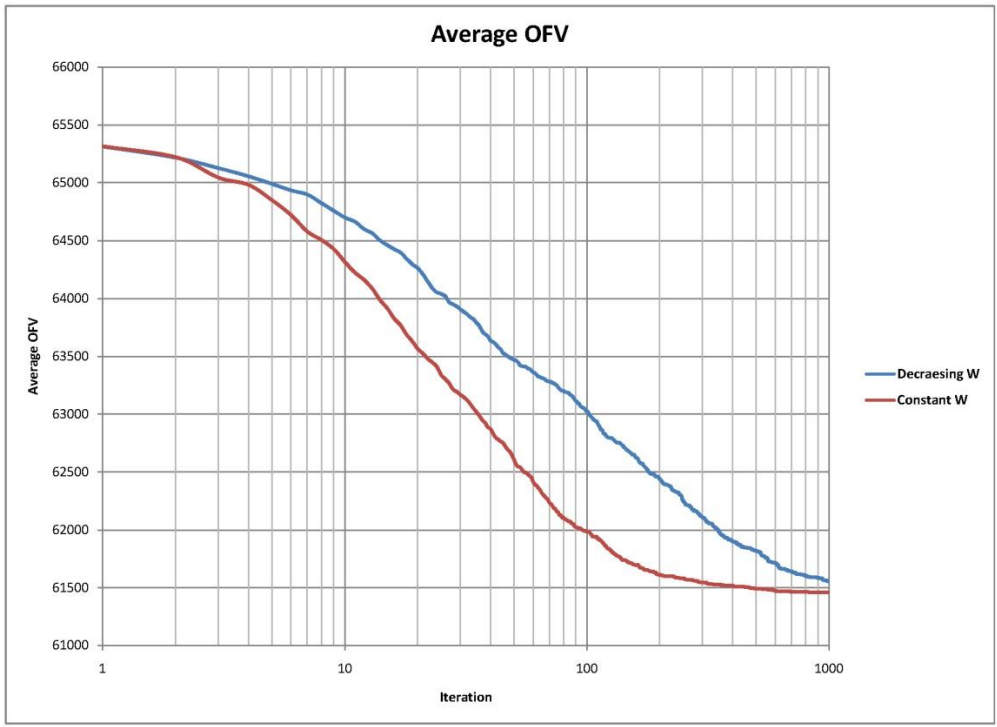

Figure 11: Average Objective Function Value for constant w and decreasing w